\documentclass[12pt]{article}                                               

\input{amssym.def}                                                           
\input{amssym.tex}                                                           
                                                                             
\begin{document}                                                             
\title{Invariants of stationary AF-algebras and  torsion
subgroup  of elliptic curves  with  complex multiplication}

\author{Igor  ~Nikolaev
\footnote{Partially supported 
by NSERC.}}


\date{}
 \maketitle


\newtheorem{thm}{Theorem}
\newtheorem{lem}{Lemma}
\newtheorem{dfn}{Definition}
\newtheorem{rmk}{Remark}
\newtheorem{cor}{Corollary}
\newtheorem{prp}{Proposition}
\newtheorem{exm}{Example}
\newtheorem{cnj}{Conjecture}

\begin{abstract}
Let $G_A$ be an  $AF$-algebra  given by 
periodic Bratteli  diagram with the incidence matrix $A\in GL(n, {\Bbb Z})$. 
For a given  polynomial $p(x)\in {\Bbb Z}[x]$  we 
assign to $G_A$ a finite abelian group $Ab_{p(x)}(G_A)={\Bbb Z}^n/p(A){\Bbb Z}^n$.
It is shown that  if  $p(0)=\pm 1$ and ${\Bbb Z}[x]/\langle p(x)\rangle$ is a principal ideal domain,
then  $Ab_{p(x)}(G_A)$ is  an invariant of the strong stable isomorphism class  
of $G_A$.  For $n=2$ and $p(x)=x-1$ 
we conjecture a formula  linking  values of the invariant  and 
 torsion  subgroup  of elliptic curves  with  complex 
multiplication.

\vspace{7mm}


\vspace{5mm}
\end{abstract}

\section{Introduction}
Let $A\in GL(n, {\Bbb Z})$ be a strictly positive  integer  
matrix and  consider the following two objects, naturally attached to $A$.
The first one, which we denote by $(G_A, \sigma_A)$, is a pair consisting 
of an $AF$-algebra, $G_A$,  given by an infinite
periodic Bratteli diagram with the incidence matrix $A$
and a shift automorphism, $\sigma_A$, canonically attached to $G_A$.
(The definitions  of an $AF$-algebra, a Bratteli diagram and a shift automorphism are given 
in Section 2.) The second object is an abelian group, which  can be introduced as
follows. Let $p(x)\in {\Bbb Z}[x]$ be a polynomial over ${\Bbb Z}$,
such that $p(0)=\pm 1$ and ${\Bbb Z}[x]/\langle p(x)\rangle$ is a principal ideal domain;
here $\langle p(x)\rangle$ means the  ideal generated by $p(x)$. 
Notice  that ${\Bbb Z}[x]/\langle p(x)\rangle$ is a principal ideal domain whenever
$p(x)$ is an irreducible polynomial and  roots of $p(x)$ generate an algebraic 
number field  whose ring of integers  is  a  principal ideal domain. 
Consider the following abelian group:
\begin{equation}\label{eq1}
{\Bbb Z}^n/p(A){\Bbb Z}^n := Ab_{p(x)}(G_A),
\end{equation}
which we shall call  an {\it abelianized} $G_A$ at the polynomial $p(x)$. 
Recall that the $AF$-algebras $G_A$ and  $G_{A'}$ are said to
be stably isomorphic,  whenever $G_A\otimes {\cal K}\cong G_{A'}\otimes {\cal K}$,
where ${\cal K}$ is the $C^*$-algebra of compact operators on a Hilbert space ${\cal H}$.
\begin{dfn}
The $AF$-algebras $G_A$ and $G_{A'}$ are said to be  strongly stably isomorphic
if they are stably isomorphic and $\sigma_A, \sigma_{A'}$ are 
the conjugate shift automorphisms. 
\end{dfn}
Roughly speaking,  the stable isomorphism is a property of $AF$-algebra $G_A$, 
while the strong stable isomorphism is a property of the $AF$-algebra $G_A$
along with  its  incidence matrix $A$. 
The main result of the present note  is the following theorem. 
\begin{thm}\label{thm1}
For each polynomial $p(x)\in {\Bbb Z}[x]$,  such that $p(0)=\pm 1$
and ${\Bbb Z}[x]/\langle p(x)\rangle$ is a principal ideal domain,
the abelian group $Ab_{p(x)}(G_A)$ is an invariant of the strong stable 
isomorphism class of the $AF$-algebra $G_A$. 
\end{thm}
\begin{rmk}
\textnormal{
The reader can find many more numerical invariants of stationary $AF$-algebras
in the remarkable monograph by [Bratteli,  Jorgensen \& Ostrovsky   2004]   \cite{BJO};
notice that  the authors consider  the case when $A$ is not necessarily a unimodular 
 matrix.  
}
\end{rmk}
Let $E_{CM}$ be an  elliptic curve with complex multiplication 
by an order of conductor $f\ge 1$ in the  imaginary quadratic field 
${\Bbb Q}(\sqrt{-d})$, where $d\ne 1$ [Silverman  1994]   \cite{S}, p. 96. 
Consider a periodic continued fraction  $f\omega=[a_0,\overline{a_1,\dots,a_n}]$,
where $\omega={1+\sqrt{d}\over 2}$  if $d\equiv 1~mod~4$ and 
$\omega=\sqrt{d}$ if $d\equiv 2,3~mod~4$.  We shall introduce 
an integer matrix $A=\prod_{i=1}^n\left(\small\matrix{a_i & 1\cr 1 & 0}\right)$,
see Section 4.1 for a motivation.  
\begin{cnj}\label{cnj1}
{\bf (``Weil's Conjecture for torsion points'')}
For each $E_{CM}$ there exists a  number field $K$
such that $E_{CM}\cong E(K)$ and a twist of $E(K)$ such that
$E_{tors}(K)\cong Ab_{x-1}(G_A)$, where $E_{tors}(K)$ is the 
torsion subgroup of $E(K)$. 
\end{cnj}
\begin{rmk}
\textnormal{
Conjecture \ref{cnj1} is an analog of (one of) classical Weil's Conjectures  for projective
varieties over finite fields,   see   e.g.  [Hartshorne 1977]  \cite{H},  pp. 449-451;
indeed,  it  identifies $E_{tors}(K)$ with the fixed points
of an automorphism $A$ of the cohomology group $H^1(E(K); {\Bbb Z})$,
see also the last paragraph of Section 3.  
}
\end{rmk}
The note is organized as follows. The preliminary facts are brought together
in Section 2.  Theorem \ref{thm1} is proved in Section 3.
In Section 4 conjecture \ref{cnj1} is  explained and some examples 
are given.

\section{Preliminaries}
 An {\it $AF$-algebra}  (approximately finite-dimensional  $C^*$-algebra) is defined to
be the  norm closure of an ascending sequence of the finite-dimensional
$C^*$-algebras $M_n$'s, where  $M_n$ is the $C^*$-algebra of the $n\times n$ matrices
with the entries in ${\Bbb C}$. Here the index $n=(n_1,\dots,n_k)$ represents
a semi-simple matrix algebra $M_n=M_{n_1}\oplus\dots\oplus M_{n_k}$.
The ascending sequence mentioned above  can be written as 
$M_1\buildrel\rm\varphi_1\over\longrightarrow M_2
   \buildrel\rm\varphi_2\over\longrightarrow\dots,
$
where $M_i$ are the finite dimensional $C^*$-algebras and
$\varphi_i$ the homomorphisms between such algebras.  The set-theoretic limit
${\cal A}=\lim M_i$ has a natural algebraic structure given by the formula
$a_m+b_k\to a+b$; here $a_m\to a,b_k\to b$ for the
sequences $a_m\in M_m,b_k\in M_k$.  
The homomorphisms $\varphi_i$ can be arranged into  a graph as follows. 
Let  $M_i=M_{i_1}\oplus\dots\oplus M_{i_k}$ and 
$M_{i'}=M_{i_1'}\oplus\dots\oplus M_{i_k'}$ be 
the semi-simple $C^*$-algebras and $\varphi_i: M_i\to M_{i'}$ the  homomorphism. 
One has the two sets of vertices $V_{i_1},\dots, V_{i_k}$ and $V_{i_1'},\dots, V_{i_k'}$
joined by the $a_{rs}$ edges, whenever the summand $M_{i_r}$ contains $a_{rs}$
copies of the summand $M_{i_s'}$ under the embedding $\varphi_i$. 
As $i$ varies, one obtains an infinite graph called a {\it Bratteli diagram} of the
$AF$-algebra   [Bratteli  1972]    \cite{Bra1}.    The Bratteli diagram defines a unique  $AF$-algebra.

If the homomorphisms $\varphi_1 =\varphi_2=\dots=Const$ in the definition of 
the $AF$-algebra ${\cal A}$,  the Bratteli diagram of  $AF$-algebra ${\cal A}$ is called {\it stationary};
by an abuse of notation,  we shall refer to the corresponding $AF$-algebra as stationary as well.   
The stationary  Bratteli diagram looks like a periodic 
graph  with the incidence matrix $A=(a_{rs})$ repeated over and over again. 
Since  matrix $A$ is a non-negative integer matrix, one can take a power of
$A$ to obtain a strictly positive integer matrix -- which we always assume 
to be the case.  We shall denote the above $AF$-algebra by $G_A$.  
Recall that in the case of $AF$-algebras, the abelian monoid $V_{\Bbb C}({\cal A})$
of finitely-generated projective modules over ${\cal A}$   (and  a scale)  
defines the $AF$-algebra up to an isomorphism and is known as a {\it dimension group}
of ${\cal A}$.  We shall use a standard dictionary existing between the $AF$-algebras and their
dimension groups [R\o rdam,  Larsen \& Laustsen  2000]   \cite{RLL},  Section 7.3.  Instead of dealing with the $AF$-algebra $G_A$,
we shall work with its dimension group $(K_0(G_A), K_0^+(G_A))$, where $K_0(G_A)$ is the lattice and 
$K_0^+(G_A)$ is a positive cone inside the lattice given by a sequence of the simplicial
dimension groups:  
\begin{equation}\label{eq2}
{\Bbb Z}^n\buildrel\rm
A
\over\longrightarrow {\Bbb Z}^n
\buildrel\rm
A
\over\longrightarrow
{\Bbb Z}^n\buildrel\rm
A
\over\longrightarrow \dots
\end{equation}
(The above notation comes from the $K_0$-group of $G_A$, see [R\o rdam,  Larsen \& Laustsen  2000]   \cite{RLL},  p.122 for the details.) 
There exists a natural automorphism, $\sigma_A$, of the dimension group
$(K_0(G_A), K_0^+(G_A))$,  see [Effros 1981]   \cite{E}, p.37. It can be defined as follows. Let $\lambda_A>1$
be the Perron-Frobenius eigenvalue and $v_A= (v_A^{(1)},\dots, v_A^{(n)})  \in {\Bbb R}^n_+$ the corresponding
eigenvector of the matrix $A$. It is known that $K_0^+(G_A)$ is defined by the
inequality ${\Bbb Z}v_A^{(1)}+\dots+{\Bbb Z}v_A^{(n)}\ge 0$ and one can
multiply ${\Bbb Z}$-module ${\Bbb Z}v_A^{(1)}+\dots+{\Bbb Z}v_A^{(n)}$
by $\lambda_A$.   It is easy to see that such a multiplication defines an automorphism
of the dimension group $(K_0(G_A), K_0^+(G_A))$. The automorphism is called a {\it shift
automorphism} and denoted by $\sigma_A$. The shift automorphisms $\sigma_A, \sigma_{A'}$
are said to be conjugate, if $\sigma_A\circ\theta=\theta\circ\sigma_{A'}$ for
some  order-isomorphism $\theta$ between the dimension groups  
$(K_0(G_A), K_0^+(G_A))$ and $(K_0(G_{A'}), K_0^+(G_{A'}))$. We shall write this fact as
$(G_A,\sigma_A)\cong (G_{A'}, \sigma_{A'})$ (an isomorphism). 
\begin{lem}\label{lm1}
The pairs $(G_A,\sigma_A)$ and $(G_{A'}, \sigma_{A'})$ are isomorphic 
if and only if the matrices $A$ and $A'$ are similar. 
\end{lem}
{\it Proof.} By Theorem 6.4 of [Effros 1981]   \cite{E}, $(G_A,\sigma_A)\cong (G_{A'}, \sigma_{A'})$
if and only if the matrices $A$ and $A'$ are shift equivalent, see [Wagoner  1999]   \cite{Wag1}
for a definition of the shift equivalence. On the other hand, since the matrices
$A$ and $A'$ are unimodular, the shift equivalence between $A$ and $A'$
coincides with a similarity of the matrices in the group $GL(n, {\Bbb Z})$,
see Corollary 2.13 of [Wagoner  1999]   \cite{Wag1}.
$\square$

\begin{cor}\label{cor1}
The $AF$-algebras $G_A$ and $G_{A'}$ are strongly stably isomorphic 
if and only if the matrices $A$ and $A'$ are similar. 
\end{cor}
{\it Proof.} By a dictionary between the dimension groups and $AF$-algebras,
the order-isomorphic dimension groups correspond to the stably isomorphic
$AF$-algebras, see Theorem 2.3. of [Effros 1981]   \cite{E}. Since $\sigma_A$ and  $\sigma_{A'}$
are conjugate, one gets a strong stable isomorphism.
$\square$

\begin{exm}
\textnormal{
Let us show that theorem \ref{thm1} is non-trivial and strong stable isomorphism 
cannot be relaxed to just stable isomorphism.  Consider the unimodular  matrices 
\begin{equation}
A=\left(\matrix{a & a-1\cr 1 & 1}\right)
~\hbox{and}
~A_h=\left(\matrix{a-h & (a-h)(h+1)-1\cr 1 & h+1}\right),
\end{equation}
where $a,h\in {\Bbb Z}$ and $a> h\ge 1$.   
Because eigenvalues of $A$ and $A_h$ coincide,  one concludes that
$(K_0(G_A), K_0^+(G_A))\cong (K_0(G_{A_h}), K_0^+(G_{A_h}))$,
i.e. $G_A$ and $G_{A_h}$ are stably isomorphic $AF$-algebras
(see Section 2 for notation). 
It is verified directly, that $\theta\circ\sigma_{A_h}=\sigma_A\circ\theta$
for $\theta=\left(\small\matrix{1 & h\cr 0 & 1}\right)$;  therefore $G_A$ and 
$G_{A_h}$ are also strongly stably isomorphic.  Notice that the strong 
stable class of $G_A$ contains more than one representative. 
Using the Smith normal form of a matrix (see below), one can find
that e.g. $Ab_{x-1}(G_A)\cong Ab_{x-1}(G_{A_h})\cong {\Bbb Z}_{a-1}$,
which is in accord with theorem \ref{thm1} for $p(x)=x-1$. 
 However,  because the eigenvalues $\lambda_A$ and $\lambda_{A^2}=\lambda_A^2$
  generate the same number field, we have an isomorphism of dimension
  groups  $(K_0(G_A), K_0^+(G_A))\cong (K_0(G_{A^2}), K_0^+(G_{A^2}))$;
  on the other hand,  because $tr~(A)\ne tr~(A^2)$  matrices $A$ and $A^2$
  (and, therefore, the shift automorphisms  $\sigma_A$ and $\sigma_{A^2}$)
  cannot be conjugate.   In this case,   the proof of theorem \ref{thm1}
 breaks,  see lemma \ref{lm1}  and Section 3; therefore strong stable isomorphism
 cannot be replaced by the stable isomorphism alone. 
}
\end{exm}

\section{Proof of theorem 1}
Our proof is based on the following criterion ([Effros 1981]   \cite{E}, Theorem 6.4):
the dimension groups 
\begin{equation}\label{eq4}
{\Bbb Z}^n\buildrel\rm
A
\over\longrightarrow {\Bbb Z}^n
\buildrel\rm
A
\over\longrightarrow
{\Bbb Z}^n\buildrel\rm
A
\over\longrightarrow \dots
\qquad
\hbox{and}
\qquad
{\Bbb Z}^n\buildrel\rm
A'
\over\longrightarrow {\Bbb Z}^n
\buildrel\rm
A'
\over\longrightarrow
{\Bbb Z}^n\buildrel\rm
A'
\over\longrightarrow \dots 
\end{equation}
are order-isomorphic and $\sigma_A, \sigma_{A'}$ are conjugate {\it iff} the matrices   
$A$ and $A'$ are similar in the group $GL(n, {\Bbb Z})$,
i.e. $A'=BAB^{-1}$ for a $B\in GL(n, {\Bbb Z})$. 
The rest of the proof follows from the structure theorem
for the finitely generated modules given by the matrix $A$ over a principal ideal domain,
 see e.g. [Shafarevich  1990]   \cite{SH}, p. 43. The result says the normal form of the module (in our case -- 
 over the principal ideal domain ${\Bbb Z}[x]/\langle p(x)\rangle$) 
is independent of the particular choice of a matrix in the similarity class of $A$.

Before proceeding to a formal proof,  let us give an intuitive idea why $Ab_{p(x)}(G_A)$
is invariant of the similarity class of  matrix $A$.   Recall that ${\Bbb Z}[x]/\langle p(x)\rangle$ 
is isomorphic to the ring of integers $O_K$ of an algebraic number field  $K={\Bbb Q}(\alpha)$,
where $\alpha$ is a root of polynomial $p(x)$.  Since $p(0)=\pm 1$ one can exclude all
rational integer entries of matrix $A\in GL(n, {\Bbb Z})$ using equation $p(\alpha)=0$;
thus one gets $A\in GL(n, O_K)$.    
But $O_K$ is a principal ideal domain (by hypothesis) and, therefore, one can use
the Euclidean algorithm to bring $A$ to a diagonal form (the Smith normal form);
the factor of $O_K$-module $GL(n, O_K)$ by a submodule defined by matrix $A$ is 
a cyclic abelian group -- denoted by $Ab_{p(x)}(G_A)$ -- which is independent 
of the similarity class of matrix $A$.   Let us pass to a step by step argument
based on the theory of modules.

By hypothesis,  ${\Bbb Z}[x]/\langle p(x)\rangle$ is a principal  ideal domain;  
we shall consider the following ${\Bbb Z}[x]/\langle p(x)\rangle$-module. If $A\in M_n({\Bbb Z})$ is an
$n\times n$ integer matrix, one endows the abelian group ${\Bbb Z}^n$ with a
${\Bbb Z}[x]/\langle p(x)\rangle$-module structure by defining:
\begin{equation}
p_n(x)v=(p_n(A))v,\quad p_n(x)\in {\Bbb Z}[x]/\langle p(x)\rangle, ~v\in {\Bbb Z}^n.
\end{equation}
 Notice that the obtained module  depends
on  matrix $A$;  we shall write  $({\Bbb Z}^n)^A$ for this module.


Fix a set of generators  $\{\varepsilon_1,\dots,\varepsilon_n\}$ of  $({\Bbb Z}^n)^A$.
We shall talk about quotient modules in terms of generators and relations,
see e.g.  lecture notes by [Morandi 2005] \cite{Mor1}.  
The relation submodule can be identified with the kernel of
a module homomorphism  $\phi_{p(x)}: ({\Bbb Z}^n)^A\to {\Bbb Z}^n$
defined  by the formula $\{p(x)\varepsilon_1,\dots, p(x)\varepsilon_n\}
\mapsto\sum_{i=1}^n p(x)\varepsilon_i$. 
The relation matrix is a mapping from the module generators 
to the relation submodule generators;  in our case the relation 
matrix is $p(A)$.   
Since relation submodule depend on the polynomial $p(x)$,
the factor-module of  ${\Bbb Z}[x]/\langle p(x)\rangle$
modulo  $\ker~\phi_{p(x)}$ will be denoted by  $({\Bbb Z}^n)^A_{p(x)}$.

Let $G=(g_{ij})$ be a matrix over the principal ideal domain,
 see e.g.  [Shafarevich  1990]   \cite{SH}, p. 43. 
 It is well known that by the elementary transformations (the Euclidean algorithm)
consisting of (i) an interchange of two rows, (ii) a multiplication of a row by $-1$,
(iii) addition a multiple of one row to another and similar operations on columns,
brings the  matrix $(g_{ij})$  to a diagonal form:
\begin{equation}
D=\left(
\matrix{
g_1 &         &     &   &       & \cr
    &  \ddots &     &   &       & \cr
    &         & g_r &   &       & \cr
    &         &     & 0 &       & \cr
    &         &     &   &\ddots & \cr
    &         &     &   &       & 0
}\right),
\end{equation}
where $g_i$ are positive integers, such that   $g_i~|~g_{i+1}$;
the latter is known the {\it Smith normal form} of a matrix over the principal ideal
domain,  see e.g.  [Shafarevich  1990]   \cite{SH}, p. 44. 
The elementary transformations are equivalent to a matrix
equation $D=PGQ$, where  $P,Q\in GL(n, {\Bbb Z})$.

We claim that  matrices $p(A)$ and $p(A')$  have the same
Smith normal form.  First,  notice that $p(A)$ and $p(A')$ are
similar matrices. Indeed, we know that
$A'$ is a matrix similar to $A$, i.e. $A'=BAB^{-1}$ for a matrix
$B\in GL(n, {\Bbb Z})$;  then it is verified directly that $p(A')=Bp(A)B^{-1}$,
i.e. $p(A)$ and $p(A')$ are similar matrices. 
Let now $D$ be the Smith normal form of $p(A)$,
then $D=Pp(A)Q$ for some $P,Q\in GL(n, {\Bbb Z})$. If $B\in  GL(n, {\Bbb Z})$
is such that $p(A')=Bp(A)B^{-1}$,  then $PB^{-1}$ and $BQ$ are also in $GL(n, {\Bbb Z})$.
One gets the following  identities:
\begin{equation}
PB^{-1}(p(A'))BQ=PB^{-1}(Bp(A)B^{-1})BQ=Pp(A)Q=D.
\end{equation}
In other words,  $p(A')$ has the same Smith normal form as $p(A)$.   
Recall that the module $({\Bbb Z}^n)^A_{p(x)}$   can be
written as:
\begin{equation}
({\Bbb Z}^n)^A_{p(x)}\cong {\Bbb Z}_{g_1}\oplus\dots\oplus {\Bbb Z}_{g_r}\oplus {\Bbb Z}^{n-r},
\end{equation}
where ${\Bbb Z}_{g_i}={\Bbb Z}/g_i {\Bbb Z}$.  Since  the same set of integers $g_i$ will appear
in the diagonal form of the matrix $p(A')$,  one gets  $Ab_{p(x)}(G_A)\cong Ab_{p(x)}(G_{A'})$
for every choice of the polynomial $p(x)$,  such that $p(0)=\pm 1$ and ${\Bbb Z}[x]/\langle p(x)\rangle$
is a principal ideal domain. 
(In the practical considerations, we often have $r=n$ so that our invariant 
is a finite abelian group.)    Theorem \ref{thm1}
follows now from  corollary \ref{cor1}.

\medskip
The most important special case of the above invariant is  when  $p(x)=x-1$
(the Bowen-Franks invariant).  The invariant takes the form:
\begin{equation}
Ab_{x-1}(G_A)={\Bbb Z}^n/(A-I){\Bbb Z}^n.
\end{equation}
The Bowen-Franks invariant is covered extensively in the literature,  see e.g. [Wagoner 1999]  \cite{Wag1};  
such an invariant  has a geometric meaning of tracking an algebraic structure of the periodic 
points of an automorphism of the lattice ${\Bbb Z}^n$ defined by the matrix $A$.
In particular, the cardinality of the group $Ab_{x-1}(G_A)$ is equal to
the total number of the isolated fixed points of the automorphism $A$.
It is easy to see, that such a number coincides with $|\det (A-I)|$.  
$\square$

\section{Torsion conjecture}
The basic facts  on elliptic curves, complex multiplication, etc.  can be found
in [Silverman  1994]   \cite{S};   an excellent introduction to the subject is 
[Silverman \& Tate  1992]  \cite{ST}.   The torsion of rational elliptic curves with complex multiplication 
was studied in [Olson  1974]   \cite{Ols1}.   A link between complex multiplication and $G_A$
was the subject of  \cite{Nik1}.   
\subsection{Teichm\"uller functor}
Let $\theta\in [0,1)$ be an irrational number.
The universal $C^*$-algebra ${\cal A}_{\theta}$ generated by the  unitaries 
$u$ and $v$ satisfying the commutation  relation $vu=e^{2\pi i\theta}uv$ 
is called a {\it noncommutative torus} [Rieffel  1990]   \cite{Rie1},
see also [Effros 1981] \cite{E}  Chapters 5  (p.34)   and  10,  and 
[R\o rdam, Larsen \& Laustsen  2000] \cite{RLL},  Exercise 5.8, pp. 86-88.  
The torus ${\cal A}_{\theta}$ is not an $AF$-algebra,   but can be embedded 
into an $AF$-algebra given by the following Bratteli diagram:

\begin{figure}[here]
\begin{picture}(300,60)(0,0)
\put(110,30){\circle{3}}
\put(120,20){\circle{3}}
\put(140,20){\circle{3}}
\put(160,20){\circle{3}}
\put(120,40){\circle{3}}
\put(140,40){\circle{3}}
\put(160,40){\circle{3}}

\put(110,30){\line(1,1){10}}
\put(110,30){\line(1,-1){10}}
\put(120,42){\line(1,0){20}}
\put(120,40){\line(1,0){20}}
\put(120,38){\line(1,0){20}}
\put(120,40){\line(1,-1){20}}
\put(120,20){\line(1,1){20}}
\put(140,41){\line(1,0){20}}
\put(140,39){\line(1,0){20}}
\put(140,40){\line(1,-1){20}}
\put(140,20){\line(1,1){20}}

\put(180,20){$\dots$}
\put(180,40){$\dots$}

\put(125,52){$a_0$}
\put(145,52){$a_1$}

\end{picture}

\caption{The $AF$-algebra  corresponding to  ${\cal A}_{\theta}$.}
\end{figure}
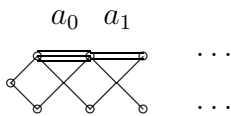

\noindent
where $\theta=[a_0, a_1,\dots]$ is the continued fraction of $\theta$,   see 
[Effros 1981] \cite{E}, p. 65.
A pair of noncommutative tori  is said to be  stably isomorphic (Morita equivalent)  whenever ${\cal A}_{\theta}\otimes {\cal K}
\cong {\cal A}_{\theta'}\otimes {\cal  K}$,  where ${\cal K}$ is the $C^*$-algebra of 
compact operators. The ${\cal A}_{\theta}$ is stably isomorphic to ${\cal A}_{\theta'}$ if
and only if  $\theta'= (a\theta +b) / (c\theta +d)$,  where  $a,b,c,d\in {\Bbb Z}$ and  $ad-bc=1$.
The K-theory of ${\cal A}_{\theta}$ is  Bott periodic with  $K_0({\cal A}_{\theta})=K_1({\cal A}_{\theta})\cong {\Bbb Z}^2$.
The range of trace on  projections of ${\cal A}_{\theta}\otimes {\cal K}$ is a subset
$\Lambda={\Bbb Z}+{\Bbb Z}\theta$  of the real line;  the set $\Lambda\cong K_0({\cal A}_{\theta})$ is known as a 
pseudo-lattice [Manin 2004]   \cite{Man1}.   The noncommutative torus ${\cal A}_{\theta}$ is said to have  {\it real multiplication},
if $\theta$ is a quadratic irrationality;  we denote such an algebra by ${\cal A}_{RM}$.
Real multiplication implies  non-trivial  endomorphisms  of the pseudo-lattice $\Lambda_{RM}$ 
given as a multiplication by real numbers -- hence the name.   Such endomorphisms 
make a ring under addition and composition of the endomorphisms;  the latter is   
isomorphic to an order of conductor $f\ge 1$ in the ring of integers of 
quadratic field ${\Bbb Q}(\theta)$.   Recall that each order of  ${\Bbb Q}(\sqrt{d})$
has the form ${\Bbb Z}+(f\omega){\Bbb Z}$,  where $\omega={1+\sqrt{d}\over 2}$ 
 if $d\equiv 1~mod~4$ and  $\omega=\sqrt{d}$ if $d\equiv 2,3~mod~4$.
It is known that continued fraction of  $\theta=f\omega$ is periodic and
has the form $[a_0,\overline{a_1,\dots,a_n}]$;  we shall consider 
a matrix $A=\prod_{i=1}^n\left(\small\matrix{a_i & 1\cr 1 & 0}\right)$. 
\begin{lem}
 $K_0(G_A)\cong K_0({\cal A}_{RM})$.
\end{lem}
{\it Proof.}
It follows easily from definition of $A$,  that $K_0(G_A)\cong {\Bbb Z}+{\Bbb Z}\theta'$,
where $\theta'=\theta-a_0$.  In other words,   $K_0(G_A)\cong K_0({\cal A}_{RM})$. 
$\square$

\bigskip\noindent
Let ${\Bbb H}=\{x+iy\in {\Bbb C}~|~y>0\}$ be the
upper half-plane and for $\tau\in {\Bbb H}$ let ${\Bbb C}/({\Bbb Z}+{\Bbb Z}\tau)$
be a complex torus;  we routinely identify the latter with a non-singular elliptic curve via the
Weierstrass $\wp$ function [Silverman  1994]   \cite{S}, pp. 6-7.  Recall that  two complex tori are 
isomorphic, whenever  $\tau'= (a\tau +b) / (c\tau +d)$,  where  $a,b,c,d\in {\Bbb Z}$ and 
$ad-bc=1$. If $\tau$ is an imaginary quadratic number,  elliptic curve is said to have 
{\it complex multiplication};  we shall denote such curves  by $E_{CM}$.
Complex multiplication means that lattice $L={\Bbb Z}+{\Bbb Z}\tau$
admits non-trivial endomorphisms given as multiplication of $L$ by certain complex (quadratic) 
numbers.  Again, such endomorphisms 
make a ring under addition and composition of the endomorphisms;  the latter is   
isomorphic to an order of conductor $f\ge 1$ in the ring of integers of 
imaginary quadratic field ${\Bbb Q}(\tau)$.

Our calculations of torsion are based on 
a covariant functor  between  elliptic curves and  noncommutative tori.
Roughly speaking, the functor maps  isomorphic curves to the stably isomorphic 
tori;  we refer the reader to  \cite{Nik1} for the details and terminology.
To give an idea,  let $\phi$ be a closed 1-form on a topological torus;
the trajectories of $\phi$  define a measured foliation on the torus.  
By the Hubbard-Masur  theorem, such a foliation 
corresponds to a point $\tau\in {\Bbb H}$.  The map $F: {\Bbb H}\to\partial {\Bbb H}$
is defined by the formula $\tau\mapsto\theta=\int_{\gamma_2}\phi/\int_{\gamma_1}\phi$,
where $\gamma_1$ and $\gamma_2$ are generators of the first homology of the
torus.    The following is true: (i) ${\Bbb H}=\partial {\Bbb H}\times (0,\infty)$
is a trivial fiber bundle, whose projection map coincides with $F$;
(ii) $F$ is a functor, which maps isomorphic complex tori to
the stably isomorphic noncommutative tori.  We shall refer to $F$
as the {\it Teichm\"uller functor}.   Remarkably,   
functor $F$  maps  $E_{CM}$  to ${\cal A}_{RM}$; 
 more specifically,  complex multiplication by order
of conductor $f$  in imaginary  field ${\Bbb Q}(\sqrt{-d})$ 
goes  to real multiplication by an order of conductor $f$
 in the real field ${\Bbb Q}(\sqrt{d})$ ,  see  an explicit
formula for $F$   \cite{Nik1}, p.524.

\subsection{Numerical examples}
We conclude by examples supporting conjecture \ref{cnj1};  they cover
all rational $E_{CM}$  [Olson  1974]   \cite{Ols1},  except $d=-1$ and $d=-163$.


\begin{table}[h]
\begin{tabular}{c|c|c|c|c|c}
\hline
        &            &    $E_{tors} ({\Bbb Q})$,   &  continued      &          &                                \\
 $-d$ &  $f$     &   see  [Olson  1974]   \cite{Ols1} &   fraction   of            & $A$  &  $Ab_{x-1}(G_A)$\\
        &            &   p.196    &   $\sqrt{f^2d}$      &          &     \\
\hline
$-2$ & $1$ & ${\Bbb Z}_2$ & $[1,\overline{2}]$ & $\left(\small\matrix{2 & 1\cr 1 & 0}\right)$ & ${\Bbb Z}_2$\\
\hline
$-3$ & $1$ & ${\Bbb Z}_1$ or ${\Bbb Z}_2$ & $[1,\overline{1,2}]$ & $\left(\small\matrix{3 & 1\cr 2 & 1}\right)$ & ${\Bbb Z}_2$\\
\hline
$-7$ & $1$ & ${\Bbb Z}_2$ & $[2,\overline{1,1,1,4}]$ & $\left(\small\matrix{14 & 3\cr 9 & 2}\right)$ & ${\Bbb Z}_{14}$\\
\hline
$-11$ & $1$  & ${\Bbb Z}_1$ & $[3,\overline{3,6}]$ & $\left(\small\matrix{19 & 3\cr 6 & 1}\right)$ & ${\Bbb Z}_3\oplus {\Bbb Z}_6$\\
\hline
$-19$ & $1$ & ${\Bbb Z}_1$ & $[4,\overline{2,1,3,1,2,8}]$ & $\left(\small\matrix{326 & 39\cr 117 & 14}\right)$ & ${\Bbb Z}_{13}\oplus {\Bbb Z}_{26}$\\
\hline
$-43$ & $1$ & ${\Bbb Z}_1$ & $[6,\overline{1,1,3,1,5,1,3,1,1,12}]$ & $\left(\small\matrix{6668 & 531\cr 3717 & 296}\right)$ & ${\Bbb Z}_{59}\oplus {\Bbb Z}_{118}$\\
\hline
$-67$ & $1$ & ${\Bbb Z}_1$ & $[8,\overline{5,2,1,1,7,1,1,2,5,16}]$ & $\left(\small\matrix{96578 & 5967\cr 17901 & 1106}\right)$ & 
${\Bbb Z}_{221}\oplus {\Bbb Z}_{442}$\\
\hline
$-3$ & $2$ & ${\Bbb Z}_2$ or ${\Bbb Z}_6$ & $[3,\overline{2,6}]$ & $\left(\small\matrix{13 & 2\cr 6 & 1}\right)$ & 
${\Bbb Z}_2\oplus {\Bbb Z}_6$\\
\hline
$-7$ & $2$ & ${\Bbb Z}_2$  & $[5,\overline{3,2,3,10}]$ & $\left(\small\matrix{247 & 24\cr 72 & 7}\right)$ & 
${\Bbb Z}_6\oplus {\Bbb Z}_{42}$\\
\hline
$-3$ & $3$ & ${\Bbb Z}_1$ & $[5,\overline{5,10}]$ & $\left(\small\matrix{51 & 5\cr 10 & 1}\right)$ & 
${\Bbb Z}_5\oplus {\Bbb Z}_{10}$\\
\hline
\end{tabular}
\end{table}
\begin{rmk}
\textnormal{
Note that $E_{tors}({\Bbb Q})\subseteq E_{tors}(K)$
since $K$ is a non-trivial extension of ${\Bbb Q}$. 
The reader can see,   that $K={\Bbb Q}$ only for the first two rows;
we do not have specific  results for $K$ in other cases,
but the table above admits existence of such a field.  
The third column lists all twists of $E({\Bbb Q})$
satisfying  conjecture \ref{cnj1}.     
}
\end{rmk}




\vskip1cm

\textsc{The Fields Institute for  Research in  Mathematical Sciences, Toronto, ON, Canada,  
E-mail:} {\sf igor.v.nikolaev@gmail.com}

\smallskip
{\it Current address: 1505-657 Worcester St.,  Southbridge,  MA 01550,  U.S.A.}

\end{document}